----------------------------------------------------------------
        \documentstyle[11pt]{article}
        \newcommand{\del}{\delta}
        \newcommand{\eps}{\epsilon}

        \newcommand{\lam}{\lambda}
        \newcommand{\om}{\omega}
        \newcommand{\th}{\theta}
        
        \newcommand{\Del}{{\it \Delta}}
        \newcommand{\Gam}{{\it \Gamma}}
        \newcommand{\PHI}{{\it \Phi}}

        \newcounter{sect}
        \newcommand{\sect}[1]{\vspace{4.5ex}\addtocounter{sect}{1}
                \begin{flushleft}
                {{\large\bf \arabic{sect}. {#1}}}
                \end{flushleft}
                \setcounter{thm}{0}
                \setcounter{equation}{0}
                \def\theequation{\arabic{sect}.\arabic{equation}}}
	\newcommand{\subsect}[2]{\vspace{2.75ex}
		\begin{flushleft}
                {{\bf \arabic{sect}.{#1} {#2}}}
                \end{flushleft}}

        \newtheorem{thm}{Theorem}[sect]
        \newtheorem{prop}[thm]{Proposition}
        \newtheorem{lemma}[thm]{Lemma}
        \newtheorem{cor}[thm]{Corollary}
        \newtheorem{defn}[thm]{Definition}
        \newtheorem{ex}[thm]{Example}
        \newtheorem{assume}[thm]{Assumption}
        \newcommand{\proof}[1]{{\em Proof.}$\quad$ {#1} $\hfill\Box$}
        \newcommand{\be}{\begin{equation}}
        \newcommand{\ee}{\end{equation}}
        \newcommand{\bea}{\begin{eqnarray}}
        \newcommand{\eea}{\end{eqnarray}}
        \newcommand{\nno}{\nonumber \\}
        \newcommand{\sep}[1]{\!\!\!\! &{#1}& \!\!\!\! }
        \newcommand{\eq}{\sep{=}}

        \newcommand{\e}{e}
        \newcommand{\I}{\sqrt{-1}\,}
        \newcommand{\bra}{\langle}
        \newcommand{\ket}{\rangle}

        \newcommand{\dr}{d}
        \newcommand{\pdr}{{\partial}}
        \newcommand{\inv}[1]{\frac{1}{#1}}
        \newcommand{\im}{{\rm Im}\,}
        \newcommand{\set}[2]{\{ {#1}\,|\, {#2}\}}
	\newcommand{\Int}[1]{{\rm Int}({#1})}
        \newcommand{\ray}{\alpha^0+t\alpha}

        \font\bbb=msym7
        \font\BBB=msym10 at 11 pt
        
        \newcommand{\RE}{{\mbox{\BBB{R}}}}
        \newcommand{\co}{{\mbox{\bbb{C}}}}
        \newcommand{\CO}{{\mbox{\BBB{C}}}}
        \font\frak=eufm10 at 12pt
        \font\sfrak=eufm9
        \newcommand{\T}{\mbox{\frak{t}}}
        \newcommand{\k}{\mbox{\frak{k}}}
        \newcommand{\g}{\mbox{\frak{g}}}
        \newcommand{\p}{\mbox{\frak{p}}}
        \newcommand{\C}{{\cal C}}
        \newcommand{\F}{{\cal F}}
        
        \newcommand{\La}{{\cal L}}
        \newcommand{\PP}{{\cal P}}
        \newcommand{\Dp}{{\cal D}'}
        \newcommand{\Sp}{{\cal S}'}

        \newcommand{\mm}{moment map }
        \newcommand{\dh}{Duistermaat-Heckman }
        \newcommand{\tdr}{\tilde{\dr}}
        \newcommand{\tom}{\tilde{\om}}
        \newcommand{\prdi}{\prod_{i=1}^n}
        \newcommand{\prdik}{\prod_{i=1}^k}
        \newcommand{\prdink}{\prod_{i=k+1}^{n}}
        \newcommand{\pap}{\prdi\alpha_i(p)}
        
        \newcommand{\papz}{{\prdi\api(\zeta)}}
        \newcommand{\papzk}{{\prdink\api(\zeta)}}
        \newcommand{\pazk}{{\prdik\alpha_i(\zeta)}}
        \newcommand{\padk}{{\prdik D_{\xi_i}}}
        \newcommand{\sump}{\sum_{p\in F}}
        \newcommand{\sumpa}{\sum_{p\in F^a}}
        \newcommand{\lvl}[2]{\frac{\om_{{#1}}^{{#2}}}{({#2})!}}
        \newcommand{\hi}[1]{H^{-1}({#1})}
        \newcommand{\poly}[1]{\bigcap_{i=1}^rK(\xi_i,{#1})}
        \newcommand{\mz}{\backslash\{0\}}
        \newcommand{\api}{\alpha^p_i}
        \newcommand{\pf}{\PHI_*|\beta|}
        \newcommand{\pfk}{J_{*}|\beta|}
        \newcommand{\feta}{{\PHI_\eta}}
        \newcommand{\fxi}{\PHI_\xi}
        \newcommand{\fxz}{\PHI_{\xi_0}}
        \newcommand{\ol}{{{\cal O}_\lam}}
        
	\newcommand{\reg}{{\rm reg}}
	\newcommand{\mma}{M\backslash M^a}
\begin{document}

$\!\,{}$

\vspace{2ex}

\begin{flushleft}
{\LARGE\bf Duistermaat-Heckman measures\\
\vspace{.5em}
in a non-compact setting}\\
\vspace{3em}
{\sc ELISA PRATO}\footnote{Address after August 1993: Ecole Normale
Sup\'erieure,
DMI, 45 Rue d'Ulm, 75230 Paris Cedex 05, France} and {\sc SIYE  WU}\\
\vspace{.5em}
{\small\em
Department of Mathematics, Princeton University, Princeton, NJ 08544, USA\\
and\\
Department of Mathematics, Columbia University, New York, NY 10027, USA}
\end{flushleft}

\vspace{2em}

\noindent
{\small {\bf Abstract.} {\em We prove a \dh type formula in a suitable
non-compact setting. We use this formula to evaluate explicitly the
pushforward of the Liouville measure via the moment map of both an
abelian and a non-abelian group action. As an application we obtain
the classical analogues of well-known multiplicity formulas for the
holomorphic discrete series representations.}}

\sect{Introduction}
Let $T$ be a torus with Lie algebra $\T$.  If $(M,\om)$ is a compact
symplectic manifold of dimension $2n$ with a Hamiltonian $T$-action,
let $\PHI\colon M\to\T^*$ be the corresponding moment map, and denote by
$\beta=\inv{(2\pi)^{2n}}\frac{\om^n}{n!}$ the Liouville volume form.
Assume initially that $M$ is compact.  Consider the integral
\be\label{Z}
\int_M\e^{\I\bra\PHI,\zeta\ket}\beta, \quad\quad\quad\zeta\in\T^\co.
	\ee
In their fundamental paper [DH] Duistermaat and Heckman use the
method of exact stationary phase to prove a formula that expresses
this integral explicitly in terms of local invariants of the $T$-fixed
point set, $F$, in $M$.

The \dh formula has a number of important applications. For example,
consider the measure $\pf$ on $\T^*$, push-forward via $\PHI$ of the
Liouville measure on $M$; we will refer to this measure as the {\em
\dh} measure. Notice that the integral (\ref{Z}) is the Fourier-Laplace
transform of $\pf$.  Guillemin, Lerman, and Sternberg [GLS] use the
\dh formula to obtain an explicit formula for $\pf$ itself under the
assumption that $F$ is isolated. This formula is generalized in [GP]
to non-abelian group actions. Recently Jeffrey and Kirwan [JK] extended
these formulas to allow non-isolated fixed points.

If $M$ is not compact the integral (\ref{Z}) may not exist.  We study
this integral in the case that there exists a component $\fxz$ of the
moment map that is proper and bounded from below.
We also assume, for simplicity, that the $T$-fixed point set is
finite; though one could more generally assume that $F$ has finitely
many connected components.  In this paper, we establish a
\dh type formula in this setting (Theorem~\ref{torus})
and as an application we obtain explicit formulas for the \dh measure
in both the abelian (Theorem~\ref{pw}) and non-abelian
(Theorem~\ref{napw}) case.  From the point of view of physics, the
integral (\ref{Z}) is the partition function of a statistical system
with phase space $M$ and energy $\fxz$. Our assumption is natural
since usually the phase space is not compact and when this is the
case, the energy is bounded from below (and not from above).

In Section 1, we explore the immediate consequences of our assumption
and review some basic facts about tempered distributions (with
typically non-compact supports) and their Fourier-Laplace
transformations.  In Section 2, we prove a \dh type formula for torus
actions; this is obtained in stages, by initially considering the case
of circle actions on manifolds with boundary.  Our formula is formally
identical to the \dh formula, except that it only makes sense for
$\im(\zeta)$ belonging to a special open cone in $\T$. (This
corresponds physically to the positivity of temperature.)  In Section
3 we first obtain an explicit formula for the measure $\pf$; then we
evaluate the measure $J_*|\beta|$, where $J$ is the moment map for the
action of a compact connected Lie group $K$ with Cartan subgroup $T$.
In Section 4, we study the regular elliptic orbits of a non-compact
semisimple Lie group $G$ that correspond to its holomorphic discrete
series representations; we observe that these orbits satisfy our
assumption and we evaluate the
\dh measures associated to the action of a compact Cartan $T$ and to
the action of a maximal compact subgroup $K$ of $G$. The non-abelian
measure was first evaluated by Duflo, Heckman and Vergne [DHV] for elliptic
orbits corresponding to all the discrete series.  Finally the appendix
contains a review of basic facts about polyhedral sets and cones that are used
throughout the paper.

	\sect{Preliminaries}
Let $M$ be a non-compact connected symplectic manifold,
$T$ a torus (with Lie algebra $\T\,$) acting on $M$ in a Hamiltonian fashion,
and $\PHI\colon M\to\T^*$ the corresponding moment map.

\subsect{1}{Some properties of the moment map}
Assume for a moment that $\fxi=\bra\PHI,\xi\ket$, for a certain $\xi\in\T$,
is a proper function (it may happen that such a $\xi$ does not exist).
Then the moment map $\PHI$ itself is a proper mapping.
Moreover, $\fxi$ is a function of Morse-Bott type with critical submanifolds
(if any) of even indices; thus the levels of $\fxi$ either are empty
or have a constant number of connected components [A].
This observation leads to strong restrictions on the occurrence of
extrema for $\fxi$, and on the image set $\fxi(M)$:

\begin{lemma}
Assume that $\fxi$ is a proper function.
If $\fxi$ is surjective there are no extrema.
If $\fxi$ is not surjective there is a unique extremal value and
$\fxi(M)$ is an interval of the types: $[m_\xi,\infty)$, $(-\infty,n_\xi]$.
\end{lemma}
\proof{According to the standard Bott-Morse theory, passing through
an extremal value entails adding an additional connected component to
the level $\fxi^{-1}(a)$; but the number of connected components is
constant so all extrema must be global.  Now, since $\fxi$ is proper
$\fxi(M)$ is an unbounded interval and there can be at most one
extremal value, none if $\fxi$ is surjective.  Assume that $\fxi$ is
not surjective.  Then we will have, for example, $\fxi>c$ strictly for
some real number $c$.  Let $a$ be a regular value of $\fxi$ and
consider the manifold $M^a=\fxi^{-1}([c,a])$ with boundary
$\fxi^{-1}(a)$; $M^a$ is compact since $\fxi$ is proper.  Let $m_\xi$
be the global minimum of $\fxi$ on $M^a$. It is easy to see that $m_\xi<a$;
this ensures that $m_\xi$ corresponds to a (global) minimum on $M$
itself, and that $\fxi(M)=[m_\xi,\infty)$.}

Let us now focus on the case where $\fxi$ is not surjective.

\begin{prop}\label{fixed}
Each connected component of the critical set of a proper component
of the moment map $\fxi$ contains at least a $T$-fixed point.
\end{prop}
\proof{Let $T'$ be the closure of the one-parameter subgroup
${\{e^{t\xi}\}}$ in $T$.  Then the critical set of $\fxi$ coincides with
$$
M_{T'}=\set{p\in M}{\mbox{stab}(p)\supseteq T'}.
$$
It is proven in [GS, Theorem~27.2], that $M_{T'}$ is
a union of connected symplectic manifolds, $Z_i$,
and that $\fxi$ maps each of these to a point.
The main observation here is that since $\fxi$ is proper $Z_i$ is compact.
But $T/T'$ acts on $Z_i$ in a Hamiltonian fashion so that
each $Z_i$ must contain points which are fixed by $T/T'$, therefore by $T$.}

Denote by $F$ the $T$-fixed point set.
{}From now on, we make the following:

\begin{assume}\label{finite}\rm
Assume that there exists a $\xi_0\in\T$ such that
$\fxz$ is proper and not surjective.
By the above proposition $F$ is non-empty; assume that it is finite.
\end{assume}

Recall that under this assumption, the image $\PHI(\mma)$ is a proper
polyhedral set\footnote{We will be using a number of properties of
{\em polyhedral sets}; we refer to the Appendix for details.}
for a sufficiently large regular value $a$ of $\fxz$ [P, Theorem~1.4].
The asymptotic cone, $\C$, of $\PHI(\mma)$ clearly does not depend on the
choice of $a$.

\begin{prop}\label{cone}
Under Assumption~\ref{finite}, let $\C\subset\T^*$ be the asymptotic cone
of the proper polyhedral set $\PHI(\mma)$, then $\fxi$ is proper
if and only if $\xi\in\pm\Int{\C'}$.	If $\xi\in\Int{\C'}$,
then $\fxi(M)= [m_\xi , \infty)$ for a suitable $m_\xi\in\RE$.
\end{prop}
\proof{Notice that $\PHI$ is proper and $\fxi=\pi_\xi\circ\PHI$,
where $\pi_\xi\colon\T^*\to\RE$ is defined by $\pi_\xi=\bra\cdot,\xi\ket$.
the proposition follows from Lemma~\ref{proj} and Lemma~\ref{bound}
since a function on $M$ is proper if and only if its restriction to
$\mma$ is.}

\subsect{2}{Distributions with non-compact support and the Laplace transform}
This subsection is devoted to a brief overview of the elements of the
theory of Laplace transforms that will be needed throughout the paper.
We refer to [H\"o] for all proofs. Let $E$ be a finite-dimensional
vector space, and let $E^*$ be its dual.  Let $\Dp(E)$ be the space of
distributions on $E$, and $\Sp(E)$, that of tempered ones.
For any $T\in\Dp(E)$,
the set $\Gam(T)=\set{\eta\in E^*}{\e^{-\bra\eta,x\ket}T(x)\in\Sp}$
(which may be empty) is convex.  Since the Fourier transform $\F$ is a
linear isomorphism from $\Sp(E)$ to $\Sp(E^*)$, if $\Gam$ is non-empty
one can define the Laplace transform of $T\in\Dp(E^*)$ by
$$
\La(T)(\zeta)=\F(\e^{-\bra\eta,x\ket}T(x))(\xi),\quad\quad\quad
\zeta=\xi+\I\eta, \quad\eta\in\Gam(T).
$$
For $T\in\Sp$, $\Gam$ contains $0$
and, by the continuity of the Fourier transform, when $\eta\to 0$
inside any closed cone in $\Gam(T)$, $\La(T)(\zeta)\to\F(T)(\xi)$ as
tempered distributions.  If the distribution is tempered and compactly
supported the region $\Gam(T)$ is all of $E$, but we will be
interested in distributions that have non-compact support.
Let's concentrate for a moment on a simple example that will be
of fundamental importance.
\begin{ex}\label{single}\rm
Let $\alpha_1,\ldots,\alpha_n$ $(n\ge\dim E)$ be a spanning set of
vectors in $E$ that generates a proper polyhedral cone, $\C_{\alpha}$.
Let $H_{\alpha_i}$ be the Heaviside distribution
defined by
$$
H_{\alpha_i}(f)=\int_0^\infty f(t\alpha_i)\dr t,
\quad\quad\quad f\in C_0^\infty(E).
$$
Then the convolution $H_{\alpha_1}\ast\cdots\ast H_{\alpha_n}$
defines a tempered distribution on $E$ supported on $\C_\alpha$.
There is another description for this measure.
Let $L_\alpha$ be the map from the positive $n$-tant in $\RE^n$ to $E$
defined by
	\be
L_\alpha(s_1,\ldots,s_n)=\sum_{i=1}^n s_i\alpha_i, \;\;\;
\mbox{where}\;\; s_i\ge 0.
	\ee
$L_\alpha$ is proper since $\C_\alpha$ is, and the
pushforward via $L_\alpha$ of Lebesgue measure, $ds$, on $E$,
is well defined and given by
$(L_\alpha)_\ast\dr s=H_{\alpha_1}\ast\cdots\ast H_{\alpha_n}$ [GP].
It is quite easy to see that the set $\Gam(L_\alpha)$ is
the dual cone $\C_\alpha'$ and that, for all $\zeta\in(E^*)^\co$ with
$\im(\zeta)\in\Int{\C_\alpha}$, the Laplace transform is given by
	\be
\La(H_{\alpha_1}\ast\cdots\ast H_{\alpha_n})(\zeta)
=\frac{(\I)^n}{\prod_{i=1}^n\bra\alpha_i, \zeta\ket}.
	\ee
\end{ex}

Returning to the Hamiltonian $T$-action on $(M^{2n},\om)$, we take
$E=\T^*$ and hence $E^*=\T$.
Let $\beta=\inv{(2\pi)^{2n}}\frac{\om^n}{n!}$ be the Liouville volume form
and $\pf$, the push-forward of the corresponding measure.
Under Assumption~\ref{finite}, $\pf$, supported on $\PHI(M)$,
is piecewise polynomial and therefore defines a tempered distribution.
It will be shown in Sect.$\,$3 that $\pf$ can be written as
a sum of distributions of the form considered in Example~\ref{single}.

\begin{prop}\label{gamma}
$\Gam(\pf)=\C'$ (the dual of the asymptotic cone, $\C$,
of $\PHI(\mma)\,$).
\end{prop}
\proof{Let $f_\PHI$ be the Radon-Nikodym derivative of $\pf$
with respect to the Lesbegue measure.
For any $\eta\in\C'$, $\e^{\I\bra\mu,\xi+\I\eta\ket}$
is bounded as $\mu$ runs through $\PHI(M)$ (Lemma~\ref{bound}).
So $\e^{\I\bra\mu,\xi+\I\eta\ket}f_\PHI(\mu)\in\Sp$, i.e., $\eta\in\Gam$.
Conversely, if $\eta\not\in\C'$, then there is an element $\alpha\in\C$
such that $\bra\alpha,\eta\ket<0$.
Moreover, from the proof of Lemma~\ref{cpp},
one can choose $\alpha^0$ such that the ray $\ray\in\PHI(\mma)$ is
contained in $\PHI(\mma)$ for sufficiently large $t>0$.
In fact, one can choose $\alpha^0$ such that the ray is in the interior
of $\PHI(\mma)$, considered as a top dimensional subset of its affine hull.
The function $f_\PHI(\ray)$ is a non-zero piecewise polynomial in $t$.
Therefore, $\e^{\I\bra\ray,\xi+\I\eta\ket}f_\PHI(\ray)$
increases at least exponentially as $t\to+\infty$.
So $\e^{\I\bra\mu,\xi+\I\eta\ket}f_\PHI(\mu)\not\in\Sp$ and $\eta\not\in\Gam$.}

	\sect{A \dh type formula}
In this section, $(M,\om)$ is a $2n$ dimensional non-compact symplectic
manifold with a Hamiltonian torus action satisfying Assumption 1.3.
For $p\in F$ let $\api$, $i=1,\ldots,n$, be the weights of the isotropy
representation of $T$ on the tangent space $T_pM$.
\begin{defn}\rm
We will say that $\zeta\in\T^{\co}$ is {\em regular} if
$$
\api(\zeta)\neq 0\;\;\;\;\mbox{for}\;\;\;\;p\in F, \;\;\;\; i=1,\ldots,n.
$$
\end{defn}
\begin{thm}\label{torus}
Under Assumption~\ref{finite}, for each regular $\zeta\in\T^{\co}$
with $\im(\zeta)\in\Int{\C'}$ we have
\be
\int_M\e^{\I\bra\PHI,\zeta\ket}\beta=(\I)^n\sump
\frac{\e^{\I\bra\PHI(p),\zeta\ket}}{\papz}.   \label{t+}
\ee
\end{thm}
\proof{Consider the lattice $\La=\set{\eta\in\T}{\e^{2\pi\eta}=1}$ in $\T$
and notice that the set
$$
{\cal A}_0 =\set{\zeta=z\eta\in\T^\co}{z\in\CO,\im(z)>0,
\eta\in\Int{\C'}\cap\La\;\;\mbox{is regular}}
$$
is dense in the set of regular elements $\zeta\in\T^{\co}$ with
$\im(\zeta)\in\Int{\C'}$.  Therefore, by continuity, it will be enough to
prove the formula for $\zeta\in{\cal A}_0$, $\bra\PHI,\zeta\ket=z\feta$.
However, since $\eta\in\La$, $H=\feta$ is a moment map for the induced
action of $S^1=\{\e^{t\eta}\}$ on $M$.  Moreover, since $\eta$ is regular,
the critical set of $\feta$ is isolated and by Proposition~\ref{fixed} it
coincides with $F$.  We have thus reduced to the special case $T=S^1$,
which follows from (\ref{c}) in Lemma~\ref{circle} below, after
taking the limit $a\to +\infty$. In fact, for $\im(z)>0$, $\e^{\I za}$
decays exponentially as $a\to +\infty$.
When $a$ is sufficiently large, the cohomological class of $\om_a$ depends
linearly on $a$ [DH], while that of $F_a$ remains fixed, since
the topology of the bundle $\hi{a}\to M_a$ does not change as $a$ runs
through a set of regular values.
So the integral over $M_a$ is a polynomial in $a$,\footnote{for an
explicit formula of this polynomial, see [W, Theorem 5.2], which can also
be deduced by collecting the coefficients of $z^{-(k+1)}$ in (\ref{c}).}
and consequently, the second sum on the right hand side of (\ref{c})
converges to $0$ as $a\to +\infty$.}

\begin{lemma}\label{circle}
Let $(M^{2n}, \om)$ be a symplectic manifold on which there is
a Hamiltonian $S^1$-action with an isolated fixed point set $F$.
Let $\alpha_1^p, \ldots, \alpha_n^p$ be the weights of the isotropy
representation of $S^1$ on $T_pM$, $p\in F$.
Assume that the \mm $H\colon M\to \RE$ is proper and not surjective.
Let $a\in\RE$ be a regular value of $H$ and let $M_a=\hi{a}/S^1$ be
the symplectic quotient with the canonical symplectic form $\om_a$.
Choose a connection of the $V$-bundle $\hi{a}\to M_a$ with curvature
2-form $F_a$. If $H$ is bounded from above (below, respectively),
let $M^a=\set{p\in M}{H(p)\ge a}$ ($\,\set{p\in M}{H(p)\le a}$, respectively)
and $F^a=F\cap M^a$. Then for any $ z\in\CO$,
        \be
\int_{M^a}\e^{\I zH}\beta=\left(\frac{\I}{z}\right)^n
\sumpa\frac{\e^{\I zH(p)}}{\pap}
\mp\inv{(2\pi)^{n-1}}\sum_{k=0}^{n-1}\frac{\e^{\I za}}{(\I z)^{k+1}}
\int_{M_a}\lvl{a}{n-1-k}\land F_a^k.			\label{c}
        \ee
\end{lemma}
\proof{Since $a$ is a regular value of $H$, there is a small number
$\del>0$ such that $\hi{(-\del,\del)}$ is diffeomorphic to
$\hi{a}\times(-\del,\del)$ and the induced symplectic form on the latter
is, up to an exact form, $\alpha\land\dr H-(H-a)F_a+\om_a$ for one
(hence any) connection 1-form $\alpha$ of $\hi{a}\to M_a$ [DH,W].
One can find an $S^1$-invariant Riemannian metric $g$ on $M$
and can choose a connection whose horizontal spaces are induced by $g$.
Denote the vector of the $S^1$-action by $X$, let $\th=i_Xg/g(X,X)$
and $\nu=\inv{(2\pi)^n}\th\land(\tdr\th)^{-1}\land\e^{\tom}$.
Here $\tdr=\dr-\I z\,i_X$ is the equivariant derivative and
$\tom=\om+\I zH$ is the closed equivariant extension of $\om$.
Both $\th$ and $\nu$ are well-defined on $M\backslash F$ and
$\inv{(2\pi)^n}\e^{\tom}=\tdr\nu$.
For any fixed point $p\in F$, let $B_p^\eps$ be the ball centered at $p$
and of radius $\eps$. Since the top form in $\inv{(2\pi)^n}\e^{\tom}$ is
$\e^{\I zH}\beta$, Stokes' theorem implies
	\be
\int_{M^a}\e^{\I zH}\beta=\sumpa\left(\int_{B_p^\eps}\e^{\I zH}\beta-
	\int_{\pdr B_p^\eps}\nu\right)+\int_{\pdr M^a}\nu.
	\ee
A standard argument [BV,GS] shows that as $\eps\to+0$,
$$\int_{B_p^\eps}\e^{\I zH}\beta-\int_{\pdr B_p^\eps}\nu=
\left(\frac{\I}{z}\right)^n\frac{\e^{\I zH(p)}}{\pap}.$$

The boundary $\pdr M^a$ is $\hi{a}$ if $H$ is bounded from above
and is $-\hi{a}$ if $H$ is bounded from below.
One can show easily that, when restricted to $\hi{a}$, $\th=\alpha$,
$\tdr\th=F_a-\I z$ and $\tom=\om_a+\I z a$. Therefore
	\bea
\int_{\pdr M^a}\nu
\eq\pm\int_{\hi{a}}\alpha\land\e^{\om_a+\I z a}\land(F_a-\I z)^{-1}	\nno
\eq\pm \frac{1}{(2\pi)^{n-1}}\left(-\inv{\I z}\right)
   \int_{M_a}\e^{\om_a+\I z a}
   \land\left(1-\frac{F_a}{\I z}\right)^{-1}				\nno
\eq\mp \frac{1}{(2\pi)^{n-1}}\sum_{k=0}^{n-1}\frac{\e^{\I z a}}{(\I z)^{k+1}}
   \int_{M_a}\lvl{a}{n-1-k}\land F_a^k.
	\eea}

Notice that if the fixed point set $F$ is no longer finite, but has
finitely many connected components, similar arguments show that
Lemma~\ref{circle}, hence Theorem~\ref{torus}, remains valid, after
replacing the point-like contribution by an integral over the fixed
(symplectic) submanifold, and the products of weights by the
equivariant Euler class of the normal bundle.

\sect{A formula for the DH measure}
\subsect{1}{The abelian case}
Now consider the Hamiltonian $T$-action on $M$.  The hyperplanes
perpendicular to the weights $\alpha_i^p$ $(p\in F, 1,\ldots,n)$
divide the cone $\C'$ into finitely many subcones, each of which we
will call a {\em Weyl chamber}.  Any regular vector $\xi\in\C'$ sits
in the interior of such a chamber.  We fix such a $\xi$ and call the
corresponding chamber, $\C^+$, the {\em positive} Weyl chamber.
Define, for $p\in F$, $i=1,\ldots,n$:
$\beta_i^p={\rm sign}(\alpha_i^p(\xi))\alpha_i^p$.
\begin{defn}\rm
The set
$\beta_i^p, p\in F, i=1,\dots,n$
is called a
{\em renormalization} of the set of weights.
\end{defn}\rm
Let $\eps(p)=\prod_{i=1}^n{\rm sign}(\alpha_i^p(\xi))$, and let
$\del_\mu$ be the delta distribution supported at $\mu\in\T^*$.

\begin{thm}\label{pw}Under Assumption~\ref{finite}
        \be
\pf=\sump\eps(p)\;\del_{\PHI(p)}
\ast H_{\beta_1^p}\ast\cdots\ast H_{\beta_n^p}. \label{gls}
        \ee
\end{thm}
Notice that (\ref{gls}) is actually a collection of formulas,
one for each choice of a positive Weyl chamber in $\C'$
and therefore of a renormalization.

\proof{Both sides of (\ref{gls}) are tempered distributions on $\T^*$.
By Theorem~\ref{torus} and Example~\ref{single}
Laplace transformations of the two sides are equal for all $\zeta$
with $\eta=\im(\zeta)$ in the interior of the positive Weyl chamber $\C^+$.
Letting $\eta\to0$ inside $\C^+$, we conclude that the
Fourier transform of the two sides of (\ref{gls}) are equal
(as tempered distributions) [H\"o].
(\ref{gls}) follows from taking the inverse Fourier transformation.}

\subsect{2}{The non-abelian case}
First some notation.  Let $K$ be a compact connected (non-abelian)
Lie group with Lie algebra $\k$ and let $T$ be a maximal torus in $K$
with Lie algebra $\T$.
Let $\Del^+=\{\I\alpha_1,\ldots,\I\alpha_k\}$ ($\alpha_i \in \T^*$)
be a set of positive roots and let $W$ be the Weyl group of the pair
$(\k^\co,\T^\co)$.
For each $i=1,\ldots,k$ let $\xi_i\in\T$ denote the vector
dual to $\alpha_i$ with respect to the Killing form and consider
$P=\prdik \xi_i$, viewed as a polynomial in $\T^*$.
Let $\T^*_{\reg}$ be the set of elements $\mu$ of $\T^*$ such that
$P(\mu)\neq 0$ and let $\k^*_{\reg}$ be the set $K\cdot\T^*_{\reg}$ in $\k^*$.

Assume that $K$ acts on a symplectic manifold $M$ in a Hamiltonian
fashion, and denote by $J$ the corresponding moment map $M \rightarrow
\k^*$. The induced action of $T$ on $M$ is also
Hamiltonian with moment map, $\PHI$, the composition of the natural
projection $\k^*
\rightarrow \T^*$ with $J$.
\begin{assume}\label{nonab}
Assume that, as a $T$-space, $M$ satisfies Assumption~\ref{finite}.
Then the $T$-fixed point set $F$ is finite; assume in addition that
$J(p)\in\k^*_{\reg}$ for each $p\in F$.
\end{assume}
Let us examine for a moment the implications of this assumption.
The first sentence implies that ($\PHI$ thus) $J$ is a proper mapping; then the
measure $\pfk$ is well defined and uniquely determined by its
$W$-invariant restriction, $\nu$, to $\T^*$, which is defined as
follows: if $f$ is a $K$-invariant smooth compactly supported
function on $\k^*$ and $g$ is its restriction to $\T^*$, then
$$
\int_{\mbox{\sfrak{t}}^*}g\;\nu=\int_{\mbox{\sfrak{k}}^*}f\;\pfk.
$$
Consider now the symplectic cross-section $X=J^{-1}(\T^*_{\reg})$; $X$
is naturally a Hamiltonian $T$-space.  The second sentence of
Assumption~\ref{nonab} implies that $X$ shares its $T$-fixed point set
with $M$ and is therefore non-empty; moreover at each of these fixed
points the set of weights of the isotropy representation of $T$ on
$T_p M$ contains the set of weights of the same representation of $T$
on $T_p X$; the weights that are left are, up to signs, the elements
$\alpha_1, \ldots, \alpha_k$.

The rest of the section will be devoted to write down and prove, under
Assumption~\ref{nonab}, an explicit formula for the measure $\nu$,
which is analogous to a formula proven in [GP] for $M$ compact.  We
begin by recalling a result of [GP], the proof of which did not rely
on the compactness of $M$.  Let $\I D_{\xi_i}$ be differentiation with
respect to $\xi_i$.
\begin{prop}[{\rm[GP]}]\label{gp}
Let $f$ be a $K$-invariant compactly supported smooth function on
$\k^*$ and let $g$ be its restriction to $\T^*$; then, for
$\zeta\in\T^\co$, 	\be \label{gpe}
\int_{\mbox{\sfrak{t}}^*} g \; e^{\I\bra - , \zeta \ket}\;\nu =
\left(\padk\right)\left(\left(\pazk\right)
\int_{\mbox{\sfrak{k}}^*} f\; e^{\I \bra - , \zeta\ket} \pfk\right).
	\ee
\end{prop}

\begin{prop} Under Assumption~\ref{nonab} $\Gam (\nu) = \C'$
(the dual of the asymptotic cone, $\C$, of $\Phi(\mma)$) and
for each $ \zeta\in\T^\co$ with $\im (\zeta)\in\C'$ we have
	\be \label{dhv}
\La(\nu)(\zeta)=
\left(\padk\right)\left(\left(\pazk\right)\La(\pfk)(\zeta)\right).
        \ee
\end{prop}
\proof{Notice that for $\zeta\in\T^\co$,
$\La (\pfk)(\zeta)=\La (\pf) (\zeta)$.
Thus from Proposition~\ref{gp} we get that $\Gam(\nu)=\Gam(\pf)$
and this last set equals $\C'$ by Proposition~\ref{gamma};
(\ref{dhv}) follows from (\ref{gpe}).}

\begin{prop}\label{nadh}
Assume that $M$ satisfies Assumption~\ref{nonab}. For each
$\zeta\in\T^\co$ with $\im (\zeta)\in\Int{\C'}$ and
$\api (\zeta)\neq 0$, $p\in F$, $i=k+1,\ldots, n$, we have
	\be\label{lnu}
\La(\nu)(\zeta)= (\I)^n \left(\padk\right)\left(\sump \eps^p
\frac{\e^{\I\bra J(p),\zeta\ket}}{\papzk}\right).
        \ee
\end{prop}
\proof{By Theorem~\ref{torus} we have, for $\zeta\in\T^\co$ regular
and $\im (\zeta)\in\Int{\C'}$,
	\be\label{usedh}
\La (\pfk)(\zeta)=(\I)^n\sump\frac{\e^{\I\bra J(p),\zeta\ket}}{\papz}.
	\ee
Now, after possible relabelings we can assume that at each $p\in F$
$\api=\eps_i^p\alpha_i$, $i= 1,\ldots,k$, with $\eps_i^p$ either $1$ or $-1$.
Denote by $\eps^p=\prdik \eps_i^p$; then by combining (\ref{dhv}) and
(\ref{usedh}) and after the appropriate cancellations we get formula
(\ref{lnu}) for the Laplace transformation of $\nu$.}

\begin{thm}\label{napw}Under Assumption~\ref{nonab}
        \be\label{na}
\nu=P\sump\eps(p)\;\del_{J(p)}\ast H_{\beta_{k+1}^p}
\ast\cdots\ast H_{\beta_n^p}.
        \ee
\end{thm}
\proof{We proceed as in the abelian case. Consider a positive Weyl
chamber $\C^+$ containing a regular element $\xi\in\C'$, and let
$\beta_i^p$, $p\in F$, $i=k+1,\ldots,n$ be the corresponding
renormalized weights (notice that here the Weyl chambers are larger
since we have deleted a number of weights.)
Let $\eps(p) = \eps^p \prod_{i=k+1}^n{\rm sign}(\alpha_i^p(\xi))$.
Now repeat step by step the proof of Theorem~\ref{pw} using
Proposition~\ref{nadh} and obtain the explicit formula of
the measure $\nu$ in (\ref{na}).}

	\sect{An application: $T$-types and $K$-types of the holomorphic
discrete series}
We begin by reviewing certain basic facts about Hermitian symmetric spaces;
we refer to [H, K1, K2] for proofs and a more extensive treatment.

Let $(G,K)$ be an {\em irreducible Hermitian symmetric} pair: $G$ is a
non-compact, simple, connected Lie group with Lie algebra $\g$ and $K$
a maximal compact subgroup with Lie algebra $\k$; $G$ has trivial
center, $K$ is connected and its center is a circle. Let $T$
be a Cartan subgroup of $K$; in this setting $T$ is
automatically a Cartan subgroup of $G$. Let $W$ denote, as in the
previous section, the Weyl group for the pair $(\k^{\co},\T^{\co})$.
It is always possible to choose a set of positive roots, $\Del^+=
\{\I\alpha_1,\ldots,\I\alpha_n\}$, for the pair $(\g^{\co},\T^{\co})$ such that
the positive non-compact roots, $\Del^+_n=\{\I\alpha_{k+1},\ldots,
\I\alpha_n\}$, are larger than the compact ones,
$\I\alpha_1,\ldots,\I\alpha_k$. With this choice, the elements of
$\Del^+_n$ agree on vectors of the one-dimensional center of $\k$; as
a consequence we have that $\Del^+_n$ is $W$-invariant and that
$\bra\alpha,\beta\ket\geq 0$ for each $\alpha,\beta\in\Del^+_n$.
Consider now the proper open $W$-symmetric cone in $\T^*$
$$
\C_n=\set{\nu\in\T^*}{\bra \alpha_i, \nu \ket >0,\;i=k+1,\ldots,n}.
$$
An {\em elliptic orbit} is, by definition, a coadjoint orbit that
intersects $\T^*$.

Consider $\lam\in\T^*_{reg}\cap\C_n$ and let $\ol$ be the elliptic
orbit through $\lam$.
The natural action of $T$ on $\ol$ is Hamiltonian with moment map,
$\PHI\colon\ol\to\T^*$, given by the restriction to $\ol$ of
the $T$-invariant projection of $\g^*$ onto $\T^*$.
Let $\xi_0$ be the unique vector in the center of $\k$ that satisfies
$\alpha_i(\xi_0)=1$ ($i=k+1,\ldots,n$; recall that such $\alpha_i$'s
agree on the center).

\begin{prop}
$\ol$ satisfies Assumptions~\ref{finite} and \ref{nonab} with respect
to $\fxz$.
\end{prop}
\proof{It is shown in [P], Propositions~2.2 and 2.3,
that $\fxz$ is proper and not surjective.
The proposition now follows from Lemma~\ref{weyl} below.}

Let $\p$ be the orthogonal complement of $\k$ in $\g$
with respect to the Killing form;
$\g = \k \oplus \p$ is called the Cartan decomposition of $\g$.
$K$ acts naturally on $\p$ via the adjoint action and the operator
$\mbox{ad}(\xi_0)$ defines a complex structure on $\p$;
from this it follows easily that $S^1=\{\e^{t\xi_0}\}$ acts freely on $\p-0$.
The $K$-equivariant diffeomorphism of $K\times\p$ onto $G$, given by
$(k,X)\to\e^Xk$, induces a $K$-equivariant diffeomorphism:
	\be\label{diffeo}
\ol\simeq K\cdot\lam\times\p.
	\ee
We then have:
\begin{lemma}\label{weyl}
The set of points of $\ol$ that are fixed by $T$ is finite and given by
$W\cdot\lam=\set{w\cdot\lam}{w\in W}\subset\T^*_{reg}$.
\end{lemma}
\proof{The points of $W\cdot\lambda$ are $T$-fixed.
Conversely, let q be a T-fixed point.
Since T is maximal abelian in G and since $S^1$, thus $T$, acts freely on
$\p-0$, one can show quite easily using the $K$-equivariant diffeomorphism
(\ref{diffeo}), that $q\in K\cdot\lam\cap\T$.
Finally, since $\lam\in\T^*_{reg}$,
$K\cdot\lam\cap\T=W\cdot\lam\subset\T^*_{reg}$.}

Notice finally that the weights of the isotropy representation of $T$
on the tangent space of $\ol$ at $w\cdot\lam$ are given by
$w\cdot\alpha_i$, $i=1,\ldots,n$.
For simplicity we renormalize the weights with respect to the chamber
containing the element $\lambda$; we then get the following corollaries of
Theorems~\ref{pw} and \ref{napw} with $\eps (p) = \eps (w)$, the determinant
of $w$ as a linear transformation of $\T^*$.

\begin{thm}
\be
\pf=\sum_{w\in W}\eps(w)\;\del_{w\cdot\lam}\ast H_{\alpha_1}\ast\cdots\ast
H_{\alpha_n}.
\ee
\end{thm}
\begin{thm}
\be
\nu=P\sum_{w\in W}\eps(w)\;
\del_{w\cdot\lam}\ast H_{\alpha_{k+1}}\ast\cdots\ast H_{\alpha_n}.
\ee
\end{thm}
The two formulas above are related to the study of the holomorphic
discrete series representations of the group $G$. In fact the first
formula is the classical analogue of a formula of Harish-Chandra [HC]
for the $T$-multiplicities of such representations. The second
formula, on the other hand, is the classical analogue of the Blattner
formula for the multiplicity of the $K$-types; the Blattner
formula gives the multiplicities of the $K$-types of all discrete
series representations.

It is quite easy to see that for the elliptic orbits that correspond
to the non-holomorphic discrete series the $T$-moment map is not
proper. This means, for example, that the measure $\pf$, in this
setting, is not well defined; this is in perfect agreement with the
representation-theoretical fact that, in this setting, the
$T$-multiplicities are not finite. On the other hand one should remark
that the $K$-moment map remains proper, that the corresponding
pushforward measure is still well defined and that it has been
explicitly evaluated for all discrete series by Duflo, Heckman, and
Vergne [DHV].  Our symplectic-theoretical approach does not extend to
the non-holomorphic case, since in our proof we are relying on the
properness of the $T$-moment map. However we are hoping that a
variation of our approach will soon yield a non-abelian formula that
holds even in the event that the $T$-moment map is not proper.

        \setcounter{sect}{1}
        \setcounter{thm}{0}\def\theequation{\Alph{sect}.\arabic{equation}}
        \setcounter{equation}{0}\def\thethm{\Alph{sect}.\arabic{thm}}
        \vspace{4.5ex}
        \begin{flushleft}
        {{\bf \Alph{sect}. Appendix: polyhedral sets and polyhedral cones}}
        \end{flushleft}
This appendix provides a self-contained account on polyhedral sets
used in the main text; we refer to [FLB] for related matters.
Let $E$ be a finite-dimensional vector space with dual $E^*$ and denote by
$\bra\cdot,\cdot\ket$ the evaluation $E^* \times E \rightarrow \RE$.
For any $\xi\in E^*$, $c\in\RE$, let $K(\xi,c)$ be the (closed) half space
$\set{\alpha\in E}{\bra\xi,\alpha\ket\ge c}$ in $E$.

\begin{defn}\rm
A {\em polyhedral set} $\PP$ in $E$ is a finite intersection
of half spaces, i.e.,
        \be
\PP=\poly{c_i}\quad\quad\mbox{for }\xi_i\in E^*,\,c_i\in\RE.
        \ee
It is called a {\em polyhedral cone} if all $c_i=0$.
\end{defn}

\begin{defn}\label{asym}\rm
Let $\PP$ be a polyhedral set in $E$. Its {\em asymptotic cone},
denoted by $\C(\PP)$, is the set of vectors $\alpha\in E$
with the property that there exists $\alpha^0\in E$
such that $\ray\in\PP$ for sufficiently large $t>0$.
\end{defn}

\begin{lemma}\label{cpp}
If $\PP=\poly{c_i}$, then $\C(\PP)=\poly{0}$.
\end{lemma}
\proof{If $\alpha\in\C(\PP)$, then for sufficiently large $t>0$,
$\ray\in\PP$, i.e, $\bra\xi_i,\ray\ket\ge c_i$.
So $\bra\xi,\alpha\ket\ge 0$, $\forall i=1,\ldots,r$,
i.e., $\alpha\in\poly{0}$.
Conversely, if $\alpha\in\poly{0}$, choose $\alpha^0\in\PP$, then
$\bra\xi_i,\ray\ket\ge c_i$ for all $t\ge 0$, i.e., $\ray\in\PP$.}

\begin{defn}\rm
A polyhedral cone $\C$ is {\em proper}
if there exists a vector $\xi\in E^*$ such that $\bra\xi,\C\mz\ket>0$.
A polyhedral set $\PP$ is {\em proper} if $\C(\PP)$ is.
\end{defn}

\begin{lemma}
A polyhedral set $\PP$ is proper if and only if it does not contain a line.
\end{lemma}
\proof{If $\PP$ contains a line $\set{\ray}{t\in\RE}$,
then $\{\pm\alpha\}\subset\C(\PP)$.
So $\C(\PP)$, hence $\PP$, is not proper.
Conversely, if $\C(\PP)$ is not proper, then $\exists\alpha\ne 0$,
$\{\pm\alpha\}\subset\C(\PP)$, which means that there exist
$\alpha^0,\alpha^1\in E$ such that $\ray,\alpha^1-t\alpha\in\PP$
for sufficiently large $t>0$.
Since $\PP$ is convex and closed, it contains the set
$\set{s\alpha^0+(1-s)\alpha^1+t\alpha}{s\in[0,1],t\in\RE}$
and hence at least a line.}

\begin{defn}\rm
The {\em dual} of a polyhedral cone $\C$ is the set
$\C'=\set{\xi\in E^*}{\bra\xi,\C\ket\ge 0}$.
\end{defn}

It is easy to see that $\C'$ is a polyhedral cone in $E^*$ and
if $\C=\poly{0}$, then $\C'=\set{\sum_{i=1}^rs_i\xi_i}{s_i\ge 0}$.

\begin{lemma}\label{bound}
For any $\xi\in E^*$, $\xi\in\C(\PP)'$
if and only if $\bra\xi,\PP\ket$ is bounded from below.
\end{lemma}
\proof{Let $\PP=\poly{c_i}$.
If $\xi\in\C(\PP)'$, then $\xi=\sum_{i=1}^rs_i\xi_i$ for some $s_i\ge 0$.
So $\bra\xi,\PP\ket=\sum_{i=1}^rs_i\bra\xi_i,\PP\ket\ge\sum_{i=1}^rs_ic_i$.
Conversely, if $\xi\not\in\C(\PP)'$, then $\exists\alpha\in\C(\PP)$,
$\bra\xi,\alpha\ket<0$.
By Definition~\ref{asym}, there exists $\alpha^0$ such that $\ray\in\PP$
for sufficiently large $t>0$.
But $\bra\xi,\ray\ket=\bra\xi,\alpha^0\ket+t\bra\xi,\alpha\ket$ is not
bounded from below.}

\begin{cor}
For any $\xi\in E^*$, $\bra\xi,\PP\ket$ is compact if and only if
$\bra\xi,\C(\PP)\ket=0$.
\end{cor}
\proof{$\bra\xi,\C(\PP)\ket=0$ is equivalent to $\{\pm\xi\}\subset\C(\PP)'$.
Using Lemma~\ref{bound}, this is equivalent to the statement that
$\bra\xi,\PP\ket$ is bounded both from above and from below.}

\begin{cor}
A polyhedral set $\PP$ is compact if and only if $\C(\PP)=\{0\}$.
\end{cor}
\proof{$\PP$ is compact if and only if for any $\xi$,
$\bra\xi,\PP\ket$ is bounded, or equivalently, $\bra\xi,\C(\PP)\ket=0$.}

\begin{lemma}\label{proj}
For any $\xi\in E^*$, let $\pi_\xi=\bra\xi,\cdot\ket\colon E\to\RE$.
Then $\pi_\xi|_\PP$ is a proper map if and only if $\xi\in\pm\Int{\C(\PP)'}$.
\end{lemma}
\proof{The inverse image
        \be
(\pi_\xi|_\PP)^{-1}([a,b])=\PP\cap\pi_\xi^{-1}([a,b])
=\PP\cap K(\xi,a)\cap K(-\xi,-b)
        \ee
is a polyhedral set with asymptotic cone
$\C_\xi=\C(\PP)\cap K(\xi,0)\cap K(-\xi,0)$.
$\pi_\xi$ is proper if and only if $\C_\xi=\{0\}$, i.e.,
for any $\alpha\in\C(\PP)\mz$, $\bra\xi,\alpha\ket\neq 0$.
This is equivalent to $\bra\xi,\C(\PP)\mz\ket>0$ or $<0$, i.e.,
$\xi\in\pm\Int{\C(\PP)'}$.}

\vspace{7ex}

\noindent{\bf Acknowledgement}\\

We would like to thank Victor Guillemin for suggesting the study
of the elliptic orbits from a symplectic viewpoint and therefore
for providing the basic motivation for the present work.

	\newpage
        \newcommand{\bib}[1]{\item[[{#1}$\!\!\!$]]}
        \newcommand{\athr}[2]{{#1}.$\,${#2}}
        \newcommand{\au}[2]{\athr{{#1}}{{#2}},}
        \newcommand{\an}[2]{\athr{{#1}}{{#2}} and}
        \newcommand{\jr}[6]{{#1}, {\it {#2}} {#3}$\,$({#4}), {#5}-{#6}.}
        \newcommand{\jn}[6]{{#1}, {\it {#2}} {#3}$\,$({#4}), {#5}-{#6};}
        \newcommand{\pr}[3]{{#1}, {#2} ({#3})}
        \newcommand{\bk}[5]{{\it {#1}}, {#2}, {#3}, {#4}, {#5}}
        \newcommand{\cf}[8]{{\it {#1}}, {#2}, {#5},
                 {#6}, {#7}, {#8}, pp.$\,${#3}-{#4}}
        \vspace{5ex}
        \begin{flushleft}
{\bf References}
        \end{flushleft}
{\small
        \begin{itemize}

        \bib{A}
        \au{M.$\,$F}{Atiyah}
        \jr{Convexity and commuting Hamiltonians}
        {Bull. London Math. Soc.}{14}{1982}{1}{15}

        \bib{BV}
        \an{N}{Berline} \au{M}{Vergne}
	\jn{Classes caract\'eristiques \'equivariantes, formule de
		localisation en cohomologie \'equivariante}
	{Comptes Rendus Acad. Sc. Paris}{295}{1982}{539}{541}
        \jr{Z\'eros d'un champ de vecteurs et classes caract\'eristiques
        \'equivariantes}{Duke Math. J.}{50}{1983}{539}{549}

	\bib{DHV}
	\au{M}{Duflo} \an{G}{Heckman} \au{M}{Vergne}
	\jr{Projection d'orbites, formule de Kirillov et formule de Blattner}
	{Mem. Soc. Math. France}{15}{1984}{65}{128}

        \bib{DH}
        \an{J.$\,$J}{Duistermaat} \au{G.$\,$J}{Heckman}
        \jn{On the variation in the cohomology of the symplectic form of the
        reduced phase space}{Invent. Math.}{69}{1982}{259}{268}
        \jr{Addendum}{{\it ibid.}}{72}{1983}{153}{158}

        \bib{FLB}
        \au{W}{Fenchel}
        \pr{Convex cones, sets, and functions}
        {Princeton University lecture notes}{1953};\\
	\au{S.$\,$R}{Lay}
	\bk{Convex sets and their applications}
	{John Wiley \& Sons}{New York, Chichester, Brisbane}{1982}
	{Chap.$\,$8};\\
	\au{A}{Br$\!\not\,$ondsted}
	\bk{An introduction to convex polytopes}
	{Springer-Verlag}{New York, Heidelberg, Berlin}{1983}{\S 8}.

        \bib{GP}
        \an{V}{Guillemin} \au{E}{Prato}
        \jr{Heckman, Kostant, and Steinberg formulas for symplectic manifolds}
        {Adv. Math.}{82}{1990}{160}{179}

%
        \bib{GS}
        \an{V}{Guillemin} \au{S}{Sternberg}
        \bk{Symplectic techniques in physics}{Cambridge University Press}
        {Cambridge, New York, Melbourne}{1990}{\S II.27}.

        \bib{He}
        \au{S}{Helgason}
        \bk{Differential geometry, Lie groups and symmetric spaces}
        {Academic Press}{Orlando, San Diego, New York}{1978}{Chap.$\,$VIII}.

        \bib{H\"o}
        \au{L}{H\"ormander}
        \bk{The analysis of linear partial differential operators I, 2nd ed.}
        {Springer-Verlag}{Berlin, Heidelberg, New York, Tokyo}{1990}{\S 7.4}.

        \bib{HC}
        Harish-Chandra,
        \jr{Representations of semisimple Lie groups IV}
        {Amer. J. Math.}{77}{1955}{743}{777}

	\bib{JK}
	\an{L.$\,$C}{Jeffrey} \au{F.$\,$C}{Kirwan}
	\pr{Localization for nonabelian group actions}
	{preprint {\tt alg-geom/ 9307001}}{1993}.

        \bib{K1}
        \au{A}{Knapp}
        \cf{Bounded symmetric domains and holomorphic discrete series}
        {Symmetric spaces}{211}{246}
        {Eds.$\;$\an{W.$\,$M}{Boothby} \au{G.$\,$L}{Weiss}}
        {Marcel Dekker, Inc.}{New York}{1972}.

        \bib{K2}
        \au{A}{Knapp}
        \bk{Representation theory of semisimple Lie groups}
        {Princeton University Press}{Princeton}{1986}{Chap.$\,$VI}.

        \bib{P}
        \au{E}{Prato}
        \pr{Convexity properties of the moment map for certain
        non-compact manifolds}{preprint}{1992}.

        \bib{W}
        \au{S}{Wu}
        \pr{An integration formula for the square of moment maps of
        circle actions}{Columbia University mathematics preprint,
	{\tt hep-th/9212071}}{1992}.

        \end{itemize}
        \end{document}